\documentclass{llncs}

\usepackage{t1enc}
\usepackage[latin2]{inputenc}
\usepackage[english]{babel}
\usepackage{amsmath}
\usepackage{amssymb}
\usepackage{bm}
\usepackage[all]{xy}%diagram
\usepackage[mathscr]{eucal}

\usepackage[normalem]{ulem}%normalem<-> \emph,\em-eket nem def-ja at alahuzasra

\newcommand{\G}{\mathscr{G}}

\title{Separation Theorem for Independent Subspace Analysis \\ with Sufficient Conditions}
\author{Zolt\'an Szab\'o, Barnab\'as P\'oczos, and Andr\'as L{\H o}rincz}
\institute{Department of Information Systems \\
E\"{o}tv\"{o}s Lor\'{a}nd University, Budapest, Hungary \\
Research Group on Intelligent Information Systems \\ Hungarian Academy of
Sciences\\
\email{szzoli@cs.elte.hu, pbarn@cs.elte.hu, lorincz@inf.elte.hu},\\
WWW home page: \url{http://nipg.inf.elte.hu}}

\begin{document}
\maketitle
\begin{abstract}
Here, a separation theorem about Independent Subspace Analysis (ISA), a generalization of Independent Component
Analysis (ICA) is proven. According to the theorem, ISA estimation can be executed in two steps under certain
conditions. In the first step, 1-dimensional ICA estimation is executed. In the second step, optimal permutation
of the ICA elements is searched for. We present sufficient conditions for the ISA Separation Theorem. Namely, we
shall show that (i) elliptically symmetric sources, (ii) 2-dimensional sources invariant to 90$^\circ$ rotation,
among others, satisfy the conditions of the theorem.
\end{abstract}

\section{Introduction}
Independent Component Analysis (ICA) \cite{jutten91blind,comon94independent}
aims to recover linearly or non-linearly mixed independent and hidden sources.
There is a broad range of applications for ICA, such as blind source separation
and blind source deconvolution \cite{bell95information}, feature extraction
\cite{bell97independent}, denoising \cite{hyvarinen99sparse}. Particular
applications include, e.g., the analysis of financial data
\cite{kiviluoto98independent}, data from neurobiology, fMRI, EEG, and MEG (see,
e.g., \cite{makeig96independent,Vigario98independent} and references therein).
For a recent review on ICA see \cite{choi05blind}.

Original ICA algorithms are 1-dimensional in the sense that all sources are assumed to be independent real
valued stochastic variables. However, applications where not all, but only certain groups of the sources are
independent may have high relevance in practice. In this case, independent sources can be multi-dimensional. For
example, consider the generalization of the cocktail-party problem, where independent groups of people are
talking about independent topics, or that more than one group of musicians are playing at the party. The
separation task requires an extension of ICA, which can be called Independent Subspace Analysis (ISA) or,
alternatively, Multi-Dimensional Independent Component Analysis (MICA)
\cite{cardoso98multidimensional,hyvarinen00emergence}. Throughout the paper, we shall use the former
abbreviation. An important application for ISA is, e.g., the processing of EEG-fMRI data \cite{akaho99MICA}.

Efforts have been made to develop ISA algorithms
\cite{cardoso98multidimensional,akaho99MICA,vollgraf01multi,bach03finding,poczos05independent1,poczos05independent2,theis05blind}. Related
theoretical problems concern mostly the estimation of entropy or mutual information. In this context, entropy estimation by Edgeworth
expansion \cite{akaho99MICA} has been extended to more than 2 dimensions and has been used for clustering and mutual information testing
\cite{hulle05edgeworth}. $k$-nearest neighbors and geodesic spanning trees methods have been applied in \cite{poczos05independent1} and
\cite{poczos05independent2} for the ISA problem. Other recent approaches search for independent subspaces via kernel methods
\cite{bach03finding} and joint block diagonalization \cite{theis05blind}.

An important observation of previous computer studies \cite{cardoso98multidimensional,poczos05independent3} is that
general ISA solver algorithms are not more efficient, in fact, sometimes produce lower quality results than simple ICA
algorithm superimposed with searches for the optimal permutation of the components. This observation led to the present
theoretical work and to some computer studies that have been published elsewhere \cite{szabo06cross}.

This technical report is constructed as follows: In
Section~\ref{sec:ISA-model} the ISA task is described.
Section~\ref{sec:ISA-sep-theorem} contains our separation theorem
for the ISA task. Sufficient conditions for the theorem are
provided in Section~\ref{sec:suff-cond}. Conclusions are drawn in
Section~\ref{sec:conclusion}.

\section{The ISA Model}\label{sec:ISA-model}
\subsection{The ISA Equations} The generative model of mixed independent multi-dimensional sources (Independent
Subspace Analysis, ISA) is the following. We assume that there are $M$ pieces of hidden $d$-dimensional sources
(\emph{components}): $\mathbf{s}^m$ \mbox{$(m=1,\ldots,M)$}. The linear transformation
\begin{eqnarray}
\mathbf{z}&=&\mathbf{A}\mathbf{s}
\end{eqnarray}
of their concatenated form
\begin{equation}
\mathbf{s}:=\left[\mathbf{s}^1;\ldots;\mathbf{s}^M\right]
\end{equation}
is available for observation only. Here, the total dimension of the sources is \mbox{$D:=d\cdot M$} and thus,
\mbox{$\mathbf{s}\in\bbbr^D$},  \mbox{$\mathbf{A}\in\bbbr^{D\times D}$} and \mbox{$\mathbf{z}\in\bbbr^D$}. In what
follows, we shall assume that \emph{mixing matrix} $\mathbf{A}$ is invertible. The ISA task is to estimate the unknown
matrix $\mathbf{A}$ (or its inverse, the so-called \emph{separation matrix} $\mathbf{W}$) and the original sources by
means of the observations $\mathbf{z}(t)$. The special case of $d=1$ corresponds to the ICA task.

\subsection{The Whiteness Assumption and its Consequences} Given our assumption on the invertibility of matrix
$\mathbf{A}$, we can assume without any loss of generality that both the sources and the observation are \emph{white},
that is,
\begin{eqnarray}
E[\mathbf{s}]&=&\mathbf{0},E\left[\mathbf{s}\mathbf{s}^T\right]=\mathbf{I}_D,\label{eq:white1}\\
E[\mathbf{z}]&=&\mathbf{0},E\left[\mathbf{z}\mathbf{z}^T\right]=\mathbf{I}_D,\label{eq:white2}
\end{eqnarray}
where superscript $T$ denotes transposition, $\mathbf{I}_D$ is the $D$-dimensional identity matrix, $E[\cdot]$ denotes
the expectation value operator. It then follows that the mixing matrix $\mathbf{A}$ and thus the separation matrix
$\mathbf{W}=\mathbf{A}^{-1}$ are orthogonal:
\begin{equation}
\mathbf{I}_D
=E\left[\mathbf{z}\mathbf{z}^T\right]=\mathbf{A}E\left[\mathbf{s}\mathbf{s}^T\right]\mathbf{A}^T=\mathbf{A}\mathbf{I}_D\mathbf{A}^T=\mathbf{A}\mathbf{A}^T.
\end{equation}
The ambiguity of the ISA task is decreased by Eqs.~\eqref{eq:white1}--\eqref{eq:white2}: Now, sources are determined up
to permutation \emph{and} orthogonal transformation of the subspaces belonging to the $\mathbf{s}^m$ sources. For more
details on this subject, see \cite{theis04uniquness}.

\subsection{The ISA Cost Function}
The ISA task can be viewed as the minimization of mutual information between the estimated components:
\begin{equation}
\min_{\mathbf{W}\in\mathcal{O}^D} I\left(\mathbf{y}^1,\ldots,\mathbf{y}^M\right)
\end{equation}
where \mbox{$\mathbf{y}=\mathbf{W}\mathbf{z}$}, \mbox{$\mathbf{y}=\left[\mathbf{y}^1;\ldots;\mathbf{y}^M\right]$} and
$\mathcal{O}^D$ denotes the space of the $D\times D$ orthogonal matrices. This cost function $I$  is equivalent to the
minimization of the sum of d-dimensional entropies, because
\begin{eqnarray}
I\left(\mathbf{y}^1,\ldots,\mathbf{y}^M\right)&=&\sum_{m=1}^MH\left(\mathbf{y}^m\right)-H(\mathbf{y})\\
 &=&\sum_{m=1}^MH\left(\mathbf{y}^m\right)-H(\mathbf{Wz})\\
 &=&\sum_{m=1}^MH\left(\mathbf{y}^m\right)-(H(\mathbf{z})+\ln(\left|\det(\mathbf{W})\right|).
\end{eqnarray}

Here, $H$ is Shannon's (multi-dimensional) differential entropy defined with logarithm of base $e$,
$\left|\cdot\right|$ denotes absolute value, `$\det$' stands for determinant. In the second equality, the
$\mathbf{y}=\mathbf{W}\mathbf{z}$ relation was exploited, and the
\begin{equation}
H(\mathbf{Wz})=H(\mathbf{z})+\ln\left(\left|\det(\mathbf{W})\right|\right)
\end{equation}
rule describing transformation of the differential entropy \cite{cover91elements} was used. $\det(\mathbf{W})=1$
because of the orthogonality of $\mathbf{W}$, so $\ln(\left|\det(\mathbf{W})\right|)=0$. The $H(\mathbf{z})$ term of
the cost is constant in $\mathbf{W}$, therefore the ISA task is equivalent to the minimization of the cost function
\begin{equation}
J(\mathbf{W}):=\sum_{m=1}^MH\left(\mathbf{y}^m\right)\label{eq:ISA-costfunction}.
\end{equation}

\section{The ISA Separation Theorem}\label{sec:ISA-sep-theorem}
The main result of this work is that the ISA task may be accomplished in two steps under certain conditions. In the
first step ICA is executed. The second step is search for the optimal permutation of the ICA components.

First, consider the so called Entropy Power Inequality (EPI)
\begin{equation}
e^{2H\left(\sum_{i=1}^Lu_i\right)}\ge \sum_{i=1}^L e^{2H(u_i)},\label{eq:EPI}
\end{equation}
where $u_1,\ldots,u_L\,\in\bbbr$ denote continuous stochastic variables. This inequality holds for example, for
independent continuous variables \cite{cover91elements}.

Let $\left\|\cdot\right\|$ denote the Euclidean norm. That is, for $\mathbf{w}\in\bbbr^L$
\begin{equation}
\left\|\mathbf{w}\right\|^2:=\sum_{i=1}^Lw_i^2,
\end{equation}
where $w_i$ is the $i^{th}$ coordinate of vector $\mathbf{w}$. The surface of the unit sphere in $L$ dimensions shall
be denoted by $S^L$:
\begin{equation}
S^L:=\{\mathbf{w}\in\bbbr^L:\left\|\mathbf{w}\right\|=1\}.
\end{equation}
If EPI is satisfied (on $S^L$) then a further inequality holds:
\begin{lemma}\label{lem:suff}
Suppose that continuous stochastic variables $u_1,\ldots,u_L \,\in\bbbr$ satisfy the following inequality
\begin{equation}
    e^{2H\left(\sum_{i=1}^Lw_iu_i\right)}\ge \sum_{i=1}^L e^{2H(w_iu_i)}, \forall\mathbf{w}\in S^L.\label{eq:w-EPI}
\end{equation}
This inequality will be called the \emph{w-EPI} condition. Then
\begin{equation}
H\left(\sum_{i=1}^L w_iu_i\right)\ge\sum_{i=1}^Lw_i^2H\left(u_i\right), \forall\mathbf{w}\in S^L.\label{eq:suff}
\end{equation}
\end{lemma}

\begin{note}
w-EPI holds, for example, for independent variables $u_i$, because independence is not affected by multiplication with
a constant.
\end{note}

\begin{proof} Assume that $\mathbf{w}\in S^L$. Applying
$\ln$ on condition \eqref{eq:w-EPI}, and using the monotonicity of the $\ln$ function, we can see that the first inequality is valid in the
following inequality chain
\begin{equation}
2H\left(\sum_{i=1}^Lw_iu_i\right)\ge \ln\left(\sum_{i=1}^L e^{2H(w_iu_i)}\right)=\ln\left(\sum_{i=1}^Le^{2H(u_i)}\cdot
w_i^2\right)\ge\sum_{i=1}^Lw_i^2\cdot\ln\left(e^{2H(u_i)}\right)=\sum_{i=1}^Lw_i^2\cdot2H(u_i).
\end{equation}
Then,
\begin{enumerate}
    \item
        we used the relation \cite{cover91elements}:
        \begin{equation}
            H(w_iu_i)=H(u_i)+\ln\left(\left|w_i\right|\right)
        \end{equation}
        for the entropy of the transformed variable. Hence
        \begin{equation}
            e^{2H(w_iu_i)}=e^{2H(u_i)+2\ln\left(\left|w_i\right|\right)}=e^{2H(u_i)}\cdot
            e^{2\ln\left(\left|w_i\right|\right)}=e^{2H(u_i)}\cdot
            w_i^2.\label{eq:2entr-transf}
        \end{equation}
    \item
        In the second inequality, we utilized the concavity of $\ln$.\qed
\end{enumerate}
\end{proof}

Now we shall use Lemma~\ref{lem:suff} to proceed. The separation theorem will be a corollary of the following claim:
\begin{proposition}\label{prop}
Let $\mathbf{y}=\left[\mathbf{y}^1;\ldots;\mathbf{y}^M\right]=\mathbf{y}(\mathbf{W})=\mathbf{W}\mathbf{s}$, where
$\mathbf{W}\in \mathcal{O}^D$, $\mathbf{y}^m$ is the estimation of the $m^{th}$ component of the ISA task. Let $y^m_i$
be the $i^{th}$ coordinate of the $m^{th}$ component. Similarly, let $s^m_i$ stand for the $i^{th}$ coordinate of the
$m^{th}$ source. Let us assume that the $\mathbf{s}^m$ sources satisfy condition~\eqref{eq:suff}. Then
\begin{equation}\label{eq:main-prop}
\sum_{m=1}^M\sum_{i=1}^dH\left(y^m_i\right)\ge
\sum_{m=1}^M\sum_{i=1}^dH\left(s^m_i\right).
\end{equation}
\end{proposition}

\begin{proof}
Let us denote the $(i,j)^{th}$ element of matrix $\mathbf{W}$ by $W_{i,j}$. Coordinates of $\mathbf{y}$ and
$\mathbf{s}$ will be denoted by $y_i$ and $s_i$, respectively. Further, let $\G^1, \ldots, \G^M$ denote the indices of
the $1^{st}, \ldots , M^{th}$ subspaces, i.e., $\G^1:=\{1,\ldots,d\},\ldots,\G^M:=\{D-d+1,\ldots,D\}$. Now, writing the
elements of the $i^{th}$ row of matrix multiplication $\mathbf{y}=\mathbf{W}\mathbf{s}$, we have
\begin{equation}
y_i=\sum_{j\in \G^1} W_{i,j}s_j+\ldots+\sum_{j\in \G^M} W_{i,j}s_j\label{eq:y=Ws}
\end{equation}
and thus,
\begin{eqnarray}
\lefteqn{H\left(y_i\right)=}\nonumber\\
&&=H\left(\sum_{j\in \G^1} W_{i,j}s_j+\ldots+\sum_{j\in \G^M} W_{i,j}s_j\right)\label{eq:H(y=Ws)}\\
&&=H\left(\left(\sum_{l\in\G^1}W_{i,l}^2\right)^{\frac{1}{2}}\frac{\sum_{j\in\G^1}W_{i,j}s_j}{\left(\sum_{l\in\G^1}W_{i,l}^2\right)^{\frac{1}{2}}}
+ \ldots + \left(\sum_{l\in\G^M}W_{i,l}^2\right)^{\frac{1}{2}}\frac{\sum_{j\in\G^M}W_{i,j}s_j}{\left(\sum_{l\in\G^M}W_{i,l}^2\right)^{\frac{1}{2}}}\right)\label{eq:w^2-in}\\
&&\ge\left(\sum_{l\in\G^1}W_{i,l}^2\right)H\left(\frac{\sum_{j\in\G^1}W_{i,j}s_j}{\left(\sum_{l\in\G^1}W_{i,l}^2\right)^{\frac{1}{2}}}\right)
+ \ldots    + \left(\sum_{l\in\G^M}W_{i,l}^2\right) H\left(\frac{\sum_{j\in\G^M}W_{i,j}s_j}{\left(\sum_{l\in\G^M}W_{i,l}^2\right)^{\frac{1}{2}}}\right)\label{eq:Lem2-applied}\\
&&=\left(\sum_{l\in\G^1}W_{i,l}^2\right)
H\left(\sum_{j\in\G^1}\frac{W_{i,j}}{\left(\sum_{l\in\G^1}W_{i,l}^2\right)^{\frac{1}{2}}}s_j\right)
+ \ldots + \left(\sum_{l\in\G^M}W_{i,l}^2\right) H\left(\sum_{j\in\G^M}\frac{W_{i,j}}{\left(\sum_{l\in\G^M}W_{i,l}^2\right)^{\frac{1}{2}}}s_j\right)\label{eq:Lem2-applied-again-pre}\\
&&\ge\left(\sum_{l\in\G^1}W_{i,l}^2\right)
\sum_{j\in\G^1}\left(\frac{W_{i,j}}{\left(\sum_{l\in\G^1}W_{i,l}^2\right)^{\frac{1}{2}}}\right)^2H\left(s_j\right)\label{eq:Lem2-applied-again}
+ \ldots + \left(\sum_{l\in\G^M}W_{i,l}^2\right) \sum_{j\in\G^M}\left(\frac{W_{i,j}}{\left(\sum_{l\in\G^M}W_{i,l}^2\right)^{\frac{1}{2}}}\right)^2H\left(s_j\right)\\
&&=\sum_{j\in\G^1}W_{i,j}^2H\left(s_j\right)+\ldots+\sum_{j\in\G^M}W_{i,j}^2H\left(s_j\right)\label{eq:H(yi-last)}
\end{eqnarray}
The above steps can be justified as follows:
\begin{enumerate}
    \item
        \eqref{eq:H(y=Ws)}: Eq.~\eqref{eq:y=Ws} was inserted into the argument of $H$.
    \item
        \eqref{eq:w^2-in}: New terms were added for Lemma~\ref{lem:suff}.
    \item
        \eqref{eq:Lem2-applied}: Sources $\mathbf{s}^m$ are independent of each other and this independence is preserved upon
        mixing \emph{within} the subspaces, and we could also use Lemma~\ref{lem:suff}, because $\mathbf{W}$
        is an orthogonal matrix.
    \item
        \eqref{eq:Lem2-applied-again-pre}: Nominators were transferred into the $\sum_j$ terms.
    \item
        \eqref{eq:Lem2-applied-again}: Variables $\mathbf{s}^m$ satisfy condition~\eqref{eq:suff} according
        to our assumptions.
    \item
        \eqref{eq:H(yi-last)}: We simplified the expression after squaring.
\end{enumerate}
Using this inequality, summing it for $i$, exchanging the order of the sums, and making use of the orthogonality of
matrix $\mathbf{W}$, we have
\begin{eqnarray}
\sum_{i=1}^DH(y_i)&\ge&\sum_{i=1}^D\left(\sum_{j\in\G^1}W_{i,j}^2H\left(s_j\right)+\ldots+\sum_{j\in\G^M}W_{i,j}^2H\left(s_j\right)\right)\\
&=&\sum_{j\in\G^1}\left(\sum_{i=1}^DW^2_{i,j}\right)H\left(s_j\right)+\ldots+\sum_{j\in\G^M}\left(\sum_{i=1}^DW^2_{i,j}\right)H\left(s_j\right)\\
&=&\sum_{j=1}^DH(s_j).
\end{eqnarray}
\qed
\end{proof}

\begin{note}
The proof holds for subspaces with different dimensions. This is also true for the following theorem.
\end{note}

Having this proposition, now we present our main theorem.
\begin{theorem}[Separation Theorem for ISA]
Presume that the $\mathbf{s}^m$ sources of the ISA model satisfy condition~\eqref{eq:suff}, and that the ICA cost
function $J(\mathbf{W})=\sum_{m=1}^M\sum_{i=1}^dH(y^m_i)$ has a minimum $\mathbf{W}\in\mathscr{O}^D$. Then it is
sufficient to search for the minimum of the ISA task as a permutation of the solution of the ICA task. Using the
concept of separation matrices, it is sufficient to explore forms
\begin{equation}
    \mathbf{W}_{\mathrm{ISA}}=\mathbf{P}\mathbf{W}_{\mathrm{ICA}},\label{eq:Wform}
\end{equation}
where $\mathbf{P} \left(\in\bbbr^{D\times D}\right)$ is a permutation matrix to be determined.
\end{theorem}
\begin{proof}
ICA minimizes the l.h.s. of Eq.~\eqref{eq:main-prop}, that is, it minimizes
\mbox{$\sum_{m=1}^M\sum_{i=1}^dH\left(y^m_i\right)$}. The set of minima is invariant to permutations and to changes of
the signs. Also, according to Proposition~\ref{prop}, $\{s^m_i\}$, i.e., the $\mathbf{s}^m$ coordinates of the
components of the solution of the ISA task belong to the set of the minima.
\end{proof}

\section{Sufficient Conditions of the Separation Theorem}\label{sec:suff-cond}
In the separation theorem, we assumed that relation~\eqref{eq:suff} is fulfilled for the $\mathbf{s}^m$ sources. Here,
we shall provide sufficient conditions when this inequality is fulfilled.

\subsection{w-EPI}
According to Lemma~\ref{lem:suff}, if the w-EPI property [i.e., \eqref{eq:w-EPI}] holds for sources $\mathbf{s}^m$,
then inequality \eqref{eq:suff} holds, too.

\subsection{Elliptically Symmetric Sources}
A stochastic variable is elliptically symmetric, or elliptical, for short, if
its density function -- which exists under mild conditions -- is constant on
elliptic surfaces.\footnote{They are often called elliptically contoured
stochastic variables.} We shall show that \eqref{eq:suff} as well as the
stronger \eqref{eq:w-EPI} w-EPI relations are fulfilled. We need certain
definitions and some basic features to prove the above statement. Thus, below
we shall elaborate on spherical (spherically symmetric) and elliptically
symmetric stochastic variables \cite{fang90symmetric,frahm04generalized}.

\subsubsection{Basic Definitions}
\begin{definition}(Characteristic function)
The characteristic function of stochastic variable \mbox{$\mathbf{v}\in\bbbr^d$} is defined by the mapping
\begin{equation}
    \bbbr^d\ni \mathbf{t}\mapsto\varphi_{\mathbf{v}}(\mathbf{t}):=E[\exp(i\mathbf{t}^T\mathbf{v})],
\end{equation}
where $i=\sqrt{-1}$ and $\exp$ is the exponential function.
\end{definition}

Spherically symmetric variables can be introduced in different ways that, together, provide the view that we need here.

\begin{definition}[Spherically symmetric variable around $\bm{\mu}$]
A stochastic variable $\mathbf{v}\in\bbbr^d$ is called spherically symmetric around $\bm{\mu}$, if:
\begin{enumerate}
    \item
        its density function is not modified by any rotation around $\bm{\mu}$. Formally, if
        \begin{equation}
        \mathbf{v}-\bm{\mu} \stackrel{\mathrm{distr}}{=}
        \mathbf{O}\left(\mathbf{v}-\bm{\mu}\right),\quad\forall \mathbf{O}\in \mathcal{O}^d,\label{eq:ss-orthinv-def}
        \end{equation}
        where $\stackrel{\mathrm{distr}}{=}$ denotes equality in distribution.
    \item
       its characteristic function with some $\phi:[0,\infty)\rightarrow\bbbr$ assumes the following form
        \begin{equation}
            \varphi_{\mathbf{v}-\bm{\mu}}(\mathbf{t})=\phi\left(\mathbf{t}^T\mathbf{t}\right).
        \end{equation} Function $\phi$ is called the \emph{characteristic generator} of $\mathbf{v}$.
    \item
        it has the following stochastic representation
        \begin{equation}
            \mathbf{v}\stackrel{\mathrm{distr}}{=}\bm{\mu}+r\mathbf{u}^{(d)},
        \end{equation}
        where
        \begin{enumerate}
            \item
                $\mathbf{\mu}\in\bbbr^d$: is a constant vector,
            \item
                $\mathbf{u}^{(d)}$: is a stochastic variable of uniform distribution over $S^d$,
            \item
                $r$: is a non-negative scalar stochastic variable, which is independent of $\mathbf{u}^{(d)}$.
        \end{enumerate}
\end{enumerate}
\end{definition}
We shall make use of the following well-known property of spherically symmetric variables:
\begin{proposition}\label{eq:prop-sph-sym-proj}
Let $\mathbf{v}$ denote a $d$-dimensional variable, which is spherically symmetric around $\bm{\mu}$. Then the
projection of $\mathbf{v}-\bm{\mu}$ onto lines through the origin have identical univariate distribution.
\end{proposition}

Affine transforms of spherically symmetric variables take us to the concept of elliptically symmetric variables.
We shall be interested in the case, when the affine transformation is bijective. Then the following definitions
are equivalent:
\begin{definition}[Elliptically symmetric variable around $\bm{\mu}$]
A stochastic variable $\mathbf{e}\in\bbbr^d$ is called elliptically symmetric around $\bm{\mu}$, if:
\begin{enumerate}
    \item
        there exists $\bm{\mu}\in\bbbr^d$ and an invertible $\bm{\Lambda}\in\bbbr^{d\times d}$ such that
        \begin{equation}
            \mathbf{e}=\bm{\mu}+\bm{\Lambda}\mathbf{v},
        \end{equation}
        where $\mathbf{v}$ is a $d$-dimensional stochastic variable, which is spherically symmetric around $\mathbf{0}$.
        In this case, the characteristic function of $\mathbf{e}$ is
        \begin{equation}
            \varphi_{\mathbf{e}}(\mathbf{t})=\exp\left(i\mathbf{t}^T\bm{\mu}\right)\phi_{\mathbf{v}}\left(\mathbf{t}^T\bm{\Sigma}\mathbf{t}\right),
        \end{equation}
        where $\bm{\Sigma}:=\bm{\Lambda}\bm{\Lambda}^T$ and $\phi_{\mathbf{v}}$ is the characteristic function
        of $\mathbf{v}$.
    \item
        there exists vector $\bm{\mu}\in\bbbr^d$, positive definite symmetric
        matrix $\bm{\Sigma}\in\bbbr^{d\times d}$, and function $\phi:[0,\infty)\rightarrow \bbbr$ such,
        that the characteristic function of $\mathbf{e}-\bm{\mu}$ is
        \begin{equation}
            \varphi_{\mathbf{e}-\bm{\mu}}(\mathbf{t})=\phi\left(\mathbf{t}^T\bm{\Sigma}\mathbf{t}\right).
        \end{equation}
        This property will be denoted as $\mathbf{e}\thicksim E_d(\bm{\mu},\bm{\Sigma},\phi)$.
        $\phi$ will be called the \emph{characteristic generator} of variable $\mathbf{e}$.
    \item
        $\mathbf{e}$ has stochastic representation of the form
        \begin{equation}
            \mathbf{e}\stackrel{\mathrm{distr}}{=}\bm{\mu}+r\bm{\Lambda}\mathbf{u}^{(d)}
        \end{equation}
        where $\bm{\Lambda}\in\bbbr^{d\times d}$ is an invertible matrix and
        \begin{enumerate}
            \item
                $\mathbf{\mu}\in\bbbr^d$: is a constant vector,
            \item
                $\mathbf{u}^{(d)}$: stochastic variable with uniform distribution on $S^d$,
            \item
                $r$: non-negative scalar stochastic variable, which is independent from $\mathbf{u}^{(d)}$.
        \end{enumerate}
\end{enumerate}
Here: $\bm{\mu}$, $\bm{\Sigma}$, and $r$ are called the \emph{location vector}, the \emph{dispersion matrix}, and the
\emph{generating variate}, respectively.
\end{definition}

\subsubsection{Basic Properties}
Here, we list important properties of an elliptic variable  $\mathbf{e}\thicksim E_d(\bm{\mu},\bm{\Sigma},\phi)$.
\begin{enumerate}
    \item
        Density function: if $\mathbf{e}$ has a density function, then it assumes the form
        \begin{equation}
            f_{\mathbf{e}}(\mathbf{x})=\left|\bm{\Lambda}\right|^{-\frac{1}{2}}\cdot
            g\left(\left(\mathbf{x}-\bm{\mu}\right)^T\bm{\Lambda}^{-1}\left(\mathbf{x}-\bm{\mu}\right)\right),\quad
            \mathbf{x}\ne\bm{\mu}\label{eq:ellipt-pdf}
        \end{equation}
        where
        \begin{equation}
            \int_0^{\infty}\frac{\pi^{\frac{d}{2}}}{\Gamma\left(\frac{d}{2}\right)}t^{\frac{d}{2}-1}g(t)\mathrm{d}t=1\label{eq:ellipt-g-cond}
        \end{equation}
        and $g:[0,\infty)\rightarrow \bbbr$ is a non-negative function.
        Here, $\Gamma$ denotes the \emph{gamma function} defined as
        \begin{equation}
            \Gamma(a):=\int_0^\infty t^a\exp(-t)\mathrm{d}t\quad(a>0).
        \end{equation}
        One can show that condition \eqref{eq:ellipt-g-cond} on $g$ is necessary and sufficient
        for making \eqref{eq:ellipt-pdf} a density function. For the existence of the density function it is
        sufficient if variable $r$ is absolutely continuous. Then function $g$
        has an explicit form, see \cite{frahm04generalized}.
    \item
        Momenta: we consider the expectation value and the variance
        \begin{equation}
            Var[\mathbf{e}]:=E\left[\left(\mathbf{e}-E[\mathbf{e}]\right)\left(\mathbf{e}-E[\mathbf{e}]\right)^T\right]
        \end{equation}
        of variable $\mathbf{e}$. They exist iff the respective momenta of $r$ are finite.
        Then, supposing that $E\left[r^2\right]$ is finite, we have
        \begin{eqnarray}
             E[\mathbf{e}]&=&\bm{\mu}\label{eq:ellipt:E}\\
             Var[\mathbf{e}]&=&\frac{E[r^2]}{d}\bm{\Sigma}=-\phi'(0)\bm{\Sigma}.\label{eq:ellipt:Var}
        \end{eqnarray}
        In what follows, we assume that $E\left[r^2\right]$ is finite.
\end{enumerate}

\subsubsection{Elliptical Sources}
Now we are ready to claim the following theorem.
\begin{proposition}
Elliptical sources $\mathbf{s}^m$ ($m=1,\ldots,M$) with finite covariances satisfy condition \eqref{eq:suff} of the ISA
separation theorem. Further, they satisfy w-EPI (with equality).
\end{proposition}
\begin{proof}
Here, we show that the w-EPI property is fulfilled with equality. Let $\mathbf{s}^m\thicksim
E_d(\bm{\mu}^m,\bm{\Sigma}^m,\phi^m)$ ($m=1,\ldots,M$) denote elliptical sources. Let us normalize each of them as
\begin{equation}
\mathbf{y}\mapsto(\bm{\Sigma}^m)^{-\frac{1}{2}}\left(\mathbf{y}-\bm{\mu}^m\right).
\end{equation}
So, it is satisfactory to prove this proposition for spherically symmetric sources. In what follows, $\mathbf{s}^m$
denotes these spherically symmetric sources. According to \eqref{eq:ellipt:E}--\eqref{eq:ellipt:Var}, spherically
symmetric sources $\mathbf{s}^m$ have zero expectation values and up to a constant multiplier they also have identity
covariance matrices:
        \begin{eqnarray}
            E[\mathbf{s}^m]&=&\mathbf{0},\\
            Var[\mathbf{s}^m]&=&c^m\cdot \mathbf{I}_d.
        \end{eqnarray}
Note that our constraint on the ISA task, namely that covariance matrices of the $\mathbf{s}^m$ sources should be equal
to $\mathbf{I}_d$, is fulfilled up to constant multipliers.

Let $P_{\mathbf{w}}$ denote the projection to straight line with direction $\mathbf{w}\in S^d$, which crosses the
origin, i.e.,
\begin{equation}
P_{\mathbf{w}}:\bbbr^d \ni\mathbf{u}\mapsto\sum_{i=1}^d w_iu_i\in\bbbr.
\end{equation}

In particular, if $\mathbf{w}$ is chosen as the canonical basis vector $\mathbf{e}_i$ (all components are 0, except the
$i^{th}$ component, which is equal to 1), then
\begin{equation}
P_{\mathbf{e_i}}(\mathbf{u})=u_i.
\end{equation}
In this interpretation, \eqref{eq:w-EPI} and w-EPI are concerned with the entropies of the projections of the different
sources onto straight lines crossing the origin. The l.h.s. projects to $\mathbf{w}$, whereas the r.h.s. projects to
the canonical basis vectors. Let $\mathbf{u}$ denote an arbitrary source, i.e., $\mathbf{u}:=\mathbf{s}^m$. According
to Proposition~\ref{eq:prop-sph-sym-proj}, distribution of the spherical $\mathbf{u}$ is the same for all such
projections and thus its entropy is identical. That is,
\begin{eqnarray}
&&\sum_{i=1}^d w_iu_i \stackrel{\mathrm{distr}}{=} u_1
\stackrel{\mathrm{distr}}{=}\ldots
\stackrel{\mathrm{distr}}{=} u_d,\quad \forall \mathbf{w}\in S^d,\\
&&H\left(\sum_{i=1}^d w_iu_i\right) = H\left(u_1\right) =\ldots = H\left(u_d\right),\quad \forall \mathbf{w}\in
S^d.\label{eq:H-inv}
\end{eqnarray}
Thus:
\begin{itemize}
    \item
        l.h.s. of w-EPI: $e^{2H(u_1)}$.
    \item
        r.h.s. of w-EPI:
        \begin{equation}
            \sum_{i=1}^d e^{2H(w_iu_i)}=\sum_{i=1}^de^{2H(u_i)}\cdot
            w_i^2=e^{2H(u_1)}\sum_{i=1}^dw_i^2=e^{2H(u_1)}\cdot 1=e^{2H(u_1)}
        \end{equation}
        At the first step, we used identity \eqref{eq:2entr-transf} for each of the terms.
        At the second step, \eqref{eq:H-inv} was utilized. Then term $e^{H(u_1)}$ was pulled out and
        we took into account that $\mathbf{w}\in S^d$.
\end{itemize}
\qed
\end{proof}
\begin{note} We note that sources of spherically symmetric distribution have already been used in the context of ISA
in \cite{hyvarinen00emergence}. In that work, a generative model was assumed. According to the assumption, the
distribution of the norms of sample projections to the subspaces were independent. This way, the task was restricted to
spherically symmetric source distributions, which is a special case of the general ISA task.
\end{note}

\subsection{Sources Invariant to $90^{\circ}$ Rotation}
In the previous section, we have seen that the case of elliptical $\mathbf{s}^m$ sources can be reduced to the
spherical case\footnote{Non-singular affine transformation can be freely performed on the sources because of the
detailed ambiguities of the ISA task.}, and that spherical variables are invariant to orthogonal transformations [see
Eq.~\eqref{eq:ss-orthinv-def}]. For mixtures of 2-dimensional components ($d=2$), much milder condition, invariance to
$90^{\circ}$ rotation, suffices. First, we observe that:

\begin{note}\label{note:ort-id-suff} In the ISA separation theorem, it is sufficient if some orthogonal
transformation of the $\mathbf{s}^m$ sources, $\mathbf{C}^m\mathbf{s}^m$ \mbox{($\mathbf{C}^m\in\mathcal{O}^d$)}
satisfy the
 condition \eqref{eq:suff}. In this case, the $\mathbf{C}^m\mathbf{s}^m$ variables are extracted by the permutation
search after the ICA transformation. Because the ISA identification has ambiguities up to orthogonal transformation in
the respective subspaces, this is suitable. In other words, for the ISA identification the existence of an Orthonormal
Basis (ONB) for each $\mathbf{u}:=\mathbf{s}^m\in\bbbr^d$ components is sufficient, on which the
\begin{equation}
h:\bbbr^d\ni\mathbf{w}\mapsto H[\left<\mathbf{w},\mathbf{u}\right>]
\end{equation}
function takes its minimum. (Here, the $\left<\mathbf{w},\mathbf{u}\right>:=\sum_{i=1}^dw_iu_i$ stochastic variable is
the projection of $\mathbf{u}$ to the direction $\mathbf{w}$.) In this case, the  entropy inequality \eqref{eq:suff} is
met with equality on the elements of the ONB.
\end{note}
Now we present our theorem concerning to the $d=2$ case.
\begin{theorem}
Let us suppose, that the density function $f$ of stochastic variable $\mathbf{u}=(u_1,u_2)(=\mathbf{s}^m)\in\bbbr^2$
exhibits the invariance
\begin{equation}
f(u_1,u_2)=f(-u_2,u_1)=f(-u_1,-u_2)=f(u_2,-u_1)\quad\left(\forall \mathbf{u}\in\bbbr^2\right),\label{sep:suff}
\end{equation}
that is, it is invariant to $90^\circ$ rotation. If function $h(\mathbf{w})=H[\left<\mathbf{w},\mathbf{u}\right>]$ has
minimum on the set $\{\mathbf{w}\ge\mathbf{0}\}\cap S^2$, it also has minimum on an ONB. \footnote{Relation
$\mathbf{w}\ge \mathbf{0}$ concerns each coordinates.} Consequently, the ISA task can be identified by the use of the
separation theorem.
\end{theorem}
\begin{proof}
Let
\begin{equation}
\mathbf{R}:=\left[\begin{array}{cc}0&-1\\1&0\end{array}\right]
\end{equation}
denote the matrix of $90^{\circ}$ ccw rotation. Let $\mathbf{w}\in S^2$. $\left<\mathbf{w},\mathbf{u}\right>\in\bbbr$
is the projection of variable $\mathbf{u}$ onto $\mathbf{w}$. The value of the density function of the stochastic
variable $\left<\mathbf{w},\mathbf{u}\right>$ in $t\in\bbbr$ (we move $t$ in direction $\mathbf{w}$) can be calculated
by integration starting from the point $\mathbf{w}t$, in direction perpendicular to $\mathbf{w}$
\begin{equation}
        f_{y=y(\mathbf{w})=\left<\mathbf{w},\mathbf{u}\right>}(t)=\int_{\mathbf{w}^\perp}
        f (\mathbf{w}t+\mathbf{z})\mathrm{d}\mathbf{z}.\label{eq:proj-int}
\end{equation}
Using the supposed invariance of $f$ and the relation \eqref{eq:proj-int} we have
\begin{equation}
f_{y(\mathbf{w})}=f_{y(\mathbf{Rw})}=f_{y(\mathbf{R}^2\mathbf{w})}=f_{y(\mathbf{R}^3\mathbf{w})},\label{eq:f-invar}
\end{equation}
where `$=$' denotes the equality of functions. Consequently, it is enough to optimize $h$ on the set $\{\mathbf{w}\ge
\mathbf{0}\}$. Let $\mathbf{w}_{min}$ be the minimum of function $h$ on the set $S^2\cap\{\mathbf{w}\ge \mathbf{0}\}$.
According to Eq.~\eqref{eq:f-invar}, $h$ takes constant and minimal values in the
\[
\{\mathbf{w}_{min},\mathbf{R}\mathbf{w}_{min},\mathbf{R}^2\mathbf{w}_{min},\mathbf{R}^3\mathbf{w}_{min}\}
\]
points. $\{\mathbf{v}_{min},\mathbf{Rv}_{min}\}$ is a suitable ONB in Note~\ref{note:ort-id-suff}.\qed
\end{proof}

\begin{note}
A special case of the requirement \eqref{sep:suff} is invariance to permutation and sign changes, that is
\begin{equation}
     f(\pm u_1,\pm u_2)=f(\pm u_2,\pm u_1).
\end{equation}
In other words, there exists a function $g:\bbbr^2\rightarrow\bbbr$, which is symmetric in its variables and
\begin{equation}
f(\mathbf{u})=g(|u_1|,|u_2|).
\end{equation}
The domain of the theorem includes
\begin{enumerate}
    \item
        the formerly presented spherical variables,
    \item
        or more generally, variables with density function of the form
        \begin{equation}
            f(\mathbf{u})=g\left(\sum_i|u_i|^p\right)\quad(p>0).
        \end{equation}
        In the literature \emph{essentially} these variables are called \emph{$L^p$-norm sphericals} (for $p>1$). Here,
        we use the \emph{$L^p$-norm spherical} denomination in a slightly extended way, for $p>0$.
\end{enumerate}
\end{note}

\subsection{Takano's Dependency Criterion}
We have seen that the w-EPI property is sufficient for the ISA separation
theorem. In \cite{takano95inequalities}, sufficient condition is provided to
satisfy the EPI condition. The condition is based on the dependencies of the
variables and it concerns the 2-dimensional case. The constraint of $d=2$ may
be generalized to higher dimensions. We are not aware of such generalizations.

We note, however, that w-EPI requires that EPI be satisfied on the surface of the unit sphere. Thus it is satisfactory
to consider the intersection of the conditions detailed in \cite{takano95inequalities} on surface of the unit sphere.

\subsection{Summary of Sufficient Conditions} Here, we summarize the presented sufficient conditions of the ISA
separation theorem. We have proven, that the requirement described by Eq.~\eqref{eq:suff} for the $\mathbf{s}^m$
sources is sufficient for the theorem. This holds if the \eqref{eq:w-EPI} w-EPI condition is fulfilled. The stronger
w-EPI is valid for
\begin{enumerate}
    \item
        sources satisfying Takano's weak dependency criterion,
    \item
        spherical sources (with equality),
    \item
        sources invariant to $90^\circ$ rotation (for $d=2$). Specially, (i) variables invariant to permutation and sign changes, and
        (ii)$L^p$-norm spherical variables belong to this family.
\end{enumerate}
These results are summarized schematically in Table~\ref{tab:suffcond-summary}.
\begin{table}
  \centering
  \caption{Sufficient conditions for the separation theorem.}\label{tab:suffcond-summary}
  \[\xymatrixcolsep{1cm}
  \xymatrix{
  &\txt{invariance to $90^{\circ}$ rotation ($d=2$)}\ar@{=>}[ddd]|-{\text{(with = for a suitable ONB)}}\ar[dr]|-{\text{specially}}&\\
  &&\text{invariance to sign and permutation}\ar[d]|-{\text{specially}}\\
  &&\text{$L^p$-norm spherical ($p>0$)}\\
  \txt{Takano's dependency\\($d=2$)} \ar@{=>}[r] & \text{w-EPI}\ar@{=>}[d] & \txt{spherical symmetry (or elliptical)}\ar@{=>}[l]_-{\text{(with = for all $\mathbf{w}\in S^d$)}} \ar[u]|-{\text{generalization for $d=2$}}\\
  &\txt{Equation~\eqref{eq:suff}: sufficient\\ for the Separation Theorem}& }
  \]
\end{table}

\section{Conclusions}\label{sec:conclusion}
In this paper a separation theorem was presented for the Independent Subspace Analysis (ISA) problem. If the
conditions of the theorem are satisfied then the ISA task can be solved in 2 steps. The first step is concerned
with the search for 1-dimensional independent components. The second step corresponds to a combinatorial
problem, the search for the optimal permutation. We have shown that elliptically symmetric sources satisfy the
conditions of the theorem. In case of 2-dimensional sources ($d=2$) invariance to $90^\circ$ rotation, or the
Takano's dependency criterion is sufficient for the separation.

These results underline our experiences that the presented 2 step procedure for solving the ISA task may produce
higher quality subspaces than sophisticated search algorithms \cite{poczos05independent1}.

Finally we mention that the possibility of this two step procedure was first noted in \cite{cardoso98multidimensional}.

\bibliographystyle{splncs}
%\bibliography{ISA_separation_theorem_TR}

\end{document}